\documentclass{amsart}
\usepackage{style}

\title{Hyperbolicity and Volumes of Bongles}
\author[C. Adams]{Colin Adams}
\address{Department of Mathematics, Williams College, Williamstown, MA 01267}
\email{cadams@williams.edu}

\author [F. Gomez-Paz]{Francisco Gomez-Paz}
\address{Department of Mathematics, MIT,  77 Massachusetts Avenue, Cambridge, MA 02139-4307} 
\email{pjgomez@mit.edu}

\author[J. Kang]{Jiachen Kang}
\address{Department of Mathematics, University of Michigan, 530 Church St, Ann Arbor, MI 48109}
\email{jiachenk@umich.edu}

\author[L. Krause]{Lukas Krause}
\address{Department of Mathematics, University of California,
970 Evans Hall, Berkeley, CA 94720-3840}
\email{lukrau2002@gmail.com}

\author[G. Li]{Gregory Li}
\address{Department of Mathematics, Harvard University, Cambridge, MA 02138}
\email{gregoryli@college.harvard.edu}

\author[R. Li]{Reyna Li}
\address{Department of Mathematics, Williams College, Williamstown, MA 01267}
\email{rl10@williams.edu}

\author[C. Marple]{Chloe Marple}
\address{Department of Mathematics, Pomona College, Claremont, CA 91711}
\email{ckme2022@mymail.pomona.edu}

\author[Z.Tan]{Ziwei Tan}
\address{Department of Mathematics, Bryn Mawr College,906 New Gulph Rd, Bryn Mawr, PA 19010}
\email{ztan2@brynmawr.edu}

\begin{document}

\begin{abstract}
     We consider a simple but infinite class of staked links known as bongles. We provide necessary and sufficient conditions for these bongles to be hyperbolic. Then, we prove that all balanced hyperbolic $n$-bongles have the same volume and the corresponding volume is an upper bound on the volume of any hyperbolic $n$-bongle for $n$ even. Moreover, all hyperbolic $n$-bongles have volume strictly less than $5n(1.01494\dots)$.  We also include explicit volume calculations for all hyperbolic 3-bongles through 6-bongles.
\end{abstract}

\maketitle 

\section{Introduction}

\begin{figure}[htbp]
\includegraphics[scale=0.6]{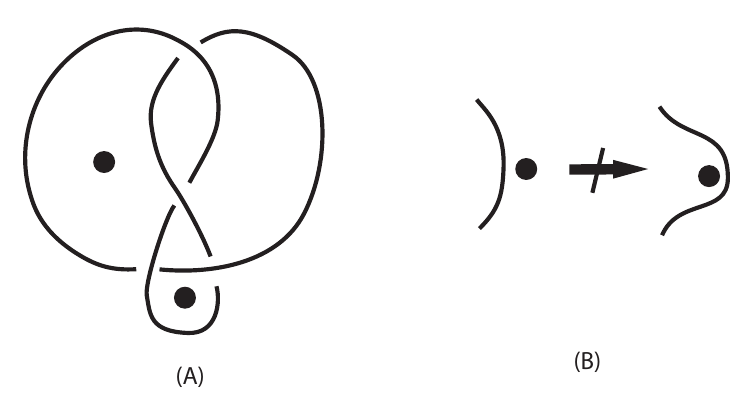}
   \caption{It is forbidden to push a strand over a stake.}
\label{fig:forbidden}
\end{figure}

A \textit{staked link} in an immersion of a collection of circles into a projection surface such that the only singularities of the immersion are transverse double points with over/under crossing data, together with a choice of points, called ``stakes," on the projection surface that avoid the immersion, as in Figure \ref{fig:forbidden}(A). The immersion is considered up to isotopy and Reidemeister moves on the surface but the strands are not allowed to pass over the stakes as in Figure \ref{fig:forbidden}(B). Staked links were introduced in full generality in \cite{generalizedknotoids} as a subset of generalized knotoids, and extend the concept of tunnel links proposed by \cite{kauffmanstaked}. 

We work with a subset of staked links on the sphere that we call \textit{charm bracelets}, which can be partitioned into \textit{charms}. In its simplest form, a charm is a part of a staked link diagram that is contained in the interior of a disk on the projection surface with boundary that is intersected exactly twice by the link diagram away from the crossings and such that the disk contains at least one crossing and at least one stake. 

As discussed in \cite{generalizedknotoids}, any staked link on the sphere can be realized as a link in a thickened sphere with open neighborhoods of the stakes removed. A thickened sphere with $n$ boundary components can then be continuously deformed to create an $(n-1)$-genus handlebody $H$. If the complement of a link in a handlebody $H$ admits a complete metric of sectional curvature $-1$ such that the handlebody boundary $\partial H$ is totally geodesic,  we say that the link complement is \textit{tg-hyperbolic}, and define the corresponding staked link to be hyperbolic.  W. Thurston proved that a link complement in a manifold is tg-hyperbolic if and only if the exterior of the link contains no properly embedded essential spheres, disks, annuli, or tori, where a sphere is essential if it does not bound a ball, and a properly embedded disk, annulus or torus is essential if it is incompressible and not boundary-parallel.

In Section 2, we introduce a bongle, which is a charm bracelet with a projection such that every charm has a single crossing that creates a monogon with a stake in its center. If there are $n$ charms, it is an $n$-bongle. 

First, we prove that a bongle is hyperbolic if and only if it is alternating. Then we consider the alternating bongles in more detail. 

Each monogon can either sit to the inside or the outside of the bongle.  In Section 3, we decompose bongles into generalized 3-bipramids and use this to obtain various results about their volumes.
First, we prove that hyperbolic $n$-bongles of fixed even $n$ with the same number of inward pointing monogons and outward pointing monogons all have identical volumes. Moreover, that volume is an upper bound on the volume of any hyperbolic $n$-bongle. And we show that 
any hyperbolic $n$-bongle has volume strictly bounded from above by $5nv_{tet}$, where $v_{tet}$ is the volume of an ideal regular hyperbolic tetrahedron, which is approximately 1.01494. We also include several conjectures about volumes of hyperbolic $n$-bongles.


For staked knots that are known to be hyperbolic, we can compute their volumes using the computer program SnapPy \cite{snappy}. We take the aforementioned handlebody $H$ containing knot $K$ and double it over its boundary to obtain a 2-component link in a connected sum of $n-1$ copies of  $S^2 \times S^1$. Such a manifold-link pair can be realized as a $(0,1)$-Dehn fillings on $n-1$ trivial components of a link of $n+1$ components in $S^3$. Taking half of the resulting volume yields the volume of the requisite manifold with totally geodesic boundary. In Section 4, we provide volumes of all hyperbolic 3-bongles through 6-bongles.

\medskip

\noindent \textbf{Acknowledgements.} \\
This research was supported by Williams College, the Finnerty Fund, and NSF Grant DMS2241623 via the SMALL Undergraduate Research Program. We are grateful to Tommy Clarke, Michael Keyes, and Felix Nusbaum, who are the other members of the SMALL 2024 knot theory group, for many helpful discussions and suggestions.

\section{Bongles}



\begin{figure}[htbp]
\includegraphics[scale=0.6]{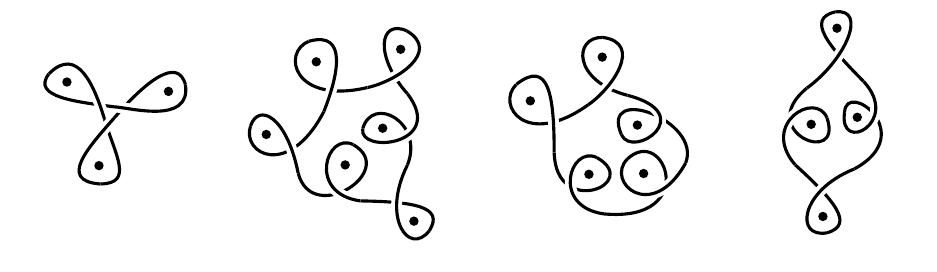}
   \caption{Examples of bongles.}
\label{fig:bonglesexamples}
\end{figure}

The simplest kinds of charm bracelets are those consisting entirely of charms, each of which has one crossing that creates a monogon containing one stake. We call these charm bracelets \textit{bongles}. Examples of bongles appear in Figure \ref{fig:bonglesexamples}. A bongle with $n$ monogons is called an $n$-bongle. Each monogon can either appear to the inside or outside of the bongle (although the choice of which region on the sphere we designate ``inside" (resp., ``outside") is arbitrary).  We call the corresponding monogons innies and outies, and a bongle with exactly the same number of innies and outies is called \textit{balanced}. When $n =2$, an $n$-bongle is equivalent to a link of two components, a case we will treat separately. 

 The first example of a hyperbolic $n$-bongle for $n \geq 3$ appeared in \cite{adamsreid}. This was an alternating 3-bongle with all monogons outies as in the first example on the left in Figure \ref{fig:bonglesexamples}. In \cite{frigerio}, this was extended to alternating $n$-bongles for $n \geq 3$ with all outies, and they were all shown to be hyperbolic. In this paper, we show that an $n$-bongle with $n \geq 2$
 is hyperbolic if and only if it is alternating. To do so in one direction, we use the following theorem from \cite{adamschen}, which provides conditions for hyperbolicity of alternating staked links. We refer to staked links and their projections interchangeably. 

We say a projection $P$ of a staked link to a surface $S$ is \textit{reduced} if there is no disk $D\subset S$ such that the interior of $D$ does not contain a stake and $\del D$ intersects the link exactly once transversely at a crossing. We say a reduced projection $P$ of a staked link is \textit{weakly prime} if there is no disk $D\subset S$ such that $\del D$ intersects the link exactly twice transversely away from the crossings and the interior of $D$ contains crossings but no stake. 

\begin{theorem}\cite{adamschen} \label{althyperbolic}
    Let $F$ be a projection surface with nonempty boundary which is not a disk, and let $L \subset F \times I$ be a link with a connected, reduced alternating projection diagram $\pi(L) \subset F \times \{1/2\}$ with at least one crossing. Let $M = (F \times I)\setminus N(L).$ Then $M$ is tg-hyperbolic if and only if the following four conditions are satisfied:
    \begin{enumerate}[label=(\roman*)]
    \item $\pi(L)$ is weakly prime on $F \times \{1/2\}$;
    \item the interior of every complementary region of $(F\times \{1/2\})\setminus \pi(L)$ is either an open disk or an open annulus;
    \item if regions $R_1$ and $R_2$ of $(F\times\{1/2\})\setminus \pi(L)$ share an edge, then at least one is a disk;
    \item there is no simple closed curve $\alpha$ in $F \times \{1/2\}$ that intersects $\pi(L)$ exactly in a nonempty collection of crossings, such that for each such crossing, $\alpha$ bisects the crossing and the two opposite complementary regions meeting at that crossing that do not intersect $\alpha$ near that crossing are annuli.
\end{enumerate}
    \label{theorem:adamschen}
\end{theorem}


We also find the following well-known result useful.

\begin{lemma} \label{boundaryparallel}If $M$ is a compact orientable irreducible $\partial$-irreducible 3-manifold, a properly embedded annulus $A$ is boundary-parallel if and only if it is $\partial$-compressible. 
\end{lemma}

\begin{proof} If $A$ is boundary-parallel, then it is clearly $\partial$-compressible. If $A$ is $\partial$-compressible, then surger the annulus along a boundary-compression disk $D$ to obtain a disk $D'$ that must also bound a disk $D''$ in the boundary of the manifold by $\partial$-irreducibility.  The two disks $D' \cup D''$ form a sphere that must bound a ball by irreducibility, which when the two copies of $D$ on its boundary are reglued, forms the solid torus that yields the fact $A$ is boundary-parallel. 
\end{proof}

We apply Theorem \ref{althyperbolic} to alternating bongles and then separately consider the non-alternating case. In our situation, $F$ is a sphere with $n$ boundary components corresponding to the $n$ stakes.
\begin{theorem}
    Let $n \geq 2$. An $n$-bongle is hyperbolic if and only if it is alternating.
\end{theorem}
\begin{proof}
We deal separately with the case of $n$ = 2. If the 2-bongle is alternating with two outies (which is equivalent to the case of two innies), the knot in the thickened sphere with two neighborhoods of stakes removed is equivalent to the Whitehead link exterior, which is well-known to be hyperbolic. If it is alternating with one outie and one innie, it is equivalent to the link exterior of the 2-component link $6_2^2$, which is also known to be hyperbolic. If the 2-bongle is not alternating, it is equivalent to either the exterior of the trivial link of two components or the exterior of the 2-braid link $6_1^2$, neither of which is hyperbolic. So, from now on, we assume $n \geq 3.$

    First, we show that all alternating bongles are hyperbolic using Theorem \ref{theorem:adamschen}. If $F$ is our projection surface and $\pi(L)$ is our projection diagram, a disk $D\subset F$ with boundary $\partial D$ that transversely intersects $\pi(L)$ exactly twice will not contain any crossings, so $\pi(L)$ is weakly prime. Condition (ii) is satisfied because each region of a bongle, by definition, can only have either zero stakes or one stake. Condition (iii) is satisfied because a bongle never has two adjacent staked regions. And Condition (iv) is satisfied because bongles never have stakes in a pair of opposite regions at a crossing. 

    Now, we show that non-alternating bongles can never be hyperbolic. We first show that every non-alternating bongle can be drawn as one of the three types of bongles depicted in Figure \ref{fig:bongleannuli} (a), allowing for the reversal of all crossings and switching the in-out order of the second pair of monogons in Type II.

    \begin{figure}[htbp]
\includegraphics[scale=0.5]{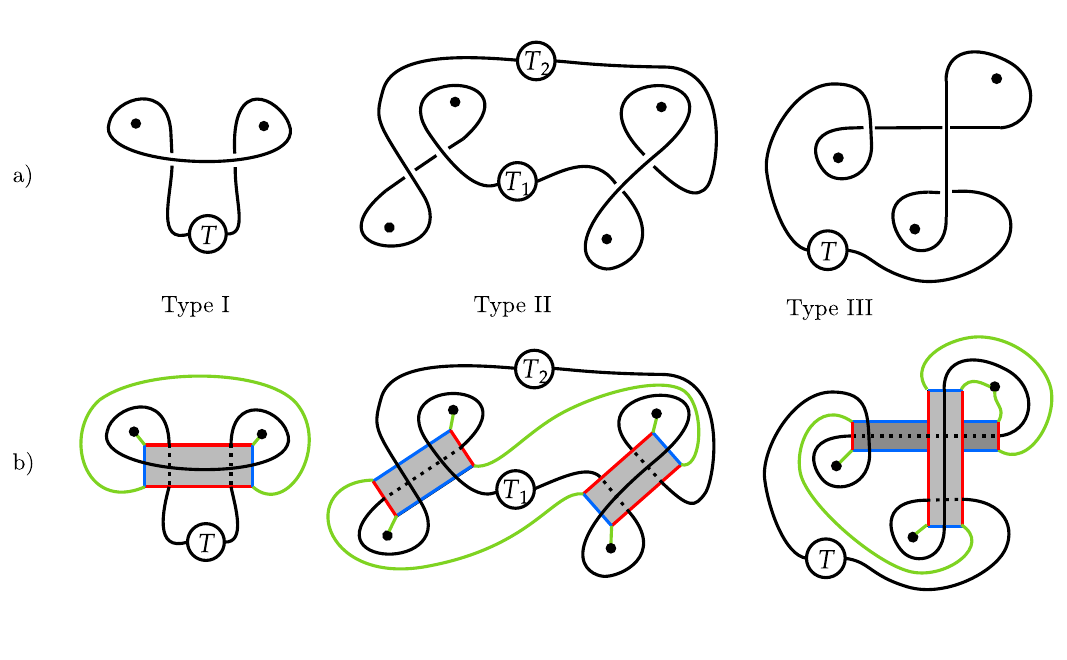}
   \caption{Essential annuli in the complement of non-alternating bongles. For simplicity, parts of the bongle are contracted into tangles, denoted $T$. }
\label{fig:bongleannuli}
\end{figure}


Begin with an arbitrary point on the knot and an arbitrary orientation. We trace around the knot and obtain a string of $u$'s and $o$'s, with $u$'s corresponding to undercrossings and $o$'s corresponding to overcrossings. Because our charms only have one crossing each, every charm we pass through adds either an $ou$ or a $uo$ to our string. Because the initial point was arbitrary, we have actually defined a cyclic string corresponding to the bongle.  
    
    Without loss of generality, a non-alternating bongle must produce two consecutive $u$'s, so assume that the first four characters of our string are $ouuo$. If we continue by adding $uo$ repeatedly to the string, then at the end, considered cyclically, we will have an $oo$. Otherwise, we switch to an $ou$ before the end, and then continue from there. Either way, for each $uu$, there must be a corresponding $oo$. 

We call a pair of monogons that are connected by a non-alternating strand a non-alternating pair. If a bongle is non-alternating, we have shown that there must be a positive even number of such pairs. 

If there is a non-alternating pair such that both the corresponding monogons are outies or both innies, we are in a bongle of Type I as in Figure \ref{fig:bongleannuli}. Otherwise, all non-alternating pairs have one monogon pointing out and one monogon pointing in. Since there must be at least two such pairs, there are two possibilities. Either the monogons in the two pairs are disjoint, or there is a single monogon that appears in both pairs. These are Types II and III bongles as in Figure \ref{fig:bongleannuli}.

    
    \begin{figure}[htbp]
\includegraphics[scale=0.5]{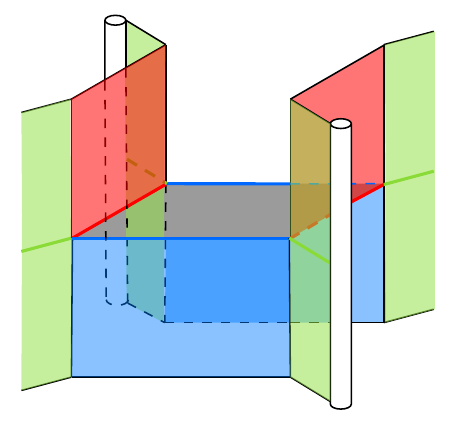}
   \caption{Example of a section of an annulus near a pair of stakes for Types II and III. For Type I, both stakes are to the front or back.}
\label{fig:bench}
\end{figure}

    Now we find an essential annulus $A$ properly embedded in the complement of each of the types of bongles shown in Figure \ref{fig:bongleannuli} (a). Consider these annuli as living locally in a thickened sphere $F\times I$ with rails removed at the stakes. In Figure \ref{fig:bongleannuli} (b), we draw the intersection of the annulus with $F\times\{1/2\}$.
    
    \begin{enumerate}
        \item Green curves shown in (b) should be interpreted as that curve cross $[0,1]$ (they extend from the top of the thickened surface to the bottom). They appear as the lightly vertical regions in Figure \ref{fig:bench}.
        \item Red curves shown in (b) should be interpreted as that edge cross $[1/2,1]$ (they extend from the middle of the thickened surface to the top). They appear as the vertical red shaded regions in Figure \ref{fig:bench}.
        \item Blue curves shown in (d) are crossed with $[0,1/2]$. They appear as the blue vertical shaded regions in Figure \ref{fig:bench}.
    \end{enumerate}
    The lightly shaded gray rectangles in Figure \ref{fig:bongleannuli} sit on the surface $F\times \{1/2\}$. They appear as the horizontal shaded regions in Figure \ref{fig:bench}. 
    
    In the case of the third bongle, there is a small neighborhood between the overpassing shaded rectangle and the underpassing one. One can think of the lighter shaded rectangle as lying in $F \times \{2/3\}$ and the darker shaded rectangle as lying in $F \times \{1/3\}$. The vertical shaded regions are adjusted accordingly.

In Figure \ref{fig:bongleannuliincomp}, we see the core curve of the corresponding annulus for each type. For Type I, the core curve is isotopic to the top boundary of the annulus $A$ and lies on the outer punctured sphere, separating the two stakes shown from the stakes that are in the tangle $T$. Hence, this is a nontrivial curve on $\partial H$ that is also nontrivial in $H$ and therefore,  $A$ is incompressible. For Types II and III, the core curve of the annulus nontrivially links the bongle link with linking number 3 in both cases and hence it must also be nontrivial, making these annuli incompressible as well.

\begin{figure}[htbp]
\includegraphics[scale=0.4]{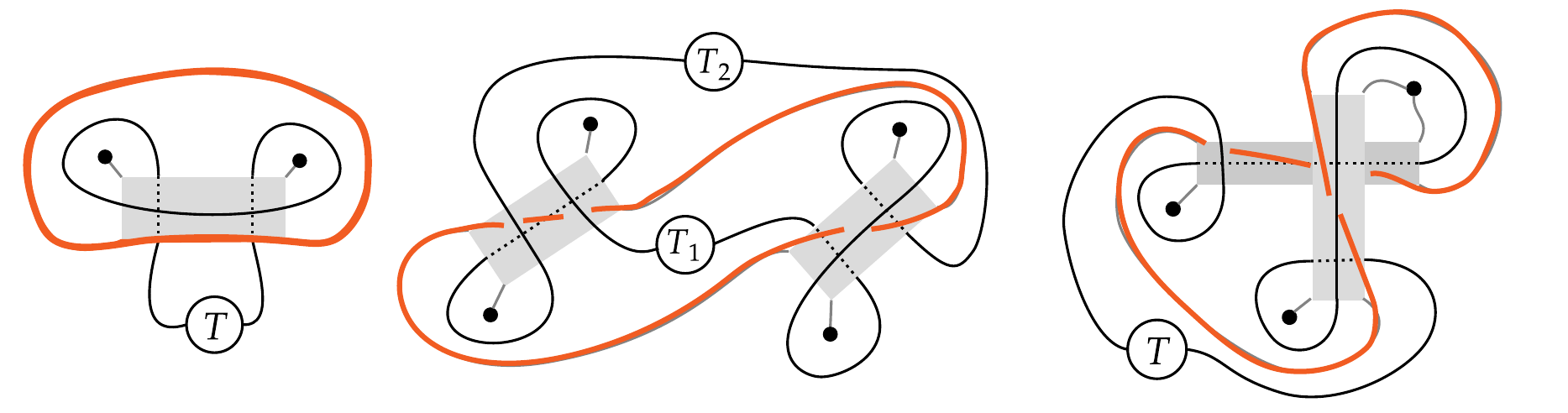}
   \caption{The core curve of each type of annulus appears in orange.}
\label{fig:bongleannuliincomp}
\end{figure}

We now prove that the bongle is not hyperbolic. Supposing the bongle is hyperbolic, the link exterior must be irreducible and $\partial$-irreducible. By Lemma \ref{boundaryparallel}, in an irreducible and $\partial$-irreducible manifold, an incompressible annulus that is $\partial$-compressible is also boundary-parallel. We see the two boundary components for each of the three types of bongles in Figure \ref{fig:bongleannulibound}, where the two blue curves form one boundary (the darker blue on $F \times \{1\}$ and the lighter blue on $F \times \{0\}$) and the two red curves form the other boundary (the darker red on $F \times \{1\}$ and the lighter red on $F \times \{0\}$). The endpoints on stakes are connected by vertical arcs on the boundary of the neighborhood of the stake we remove.


\begin{figure}[htbp]
\includegraphics[scale=0.4]{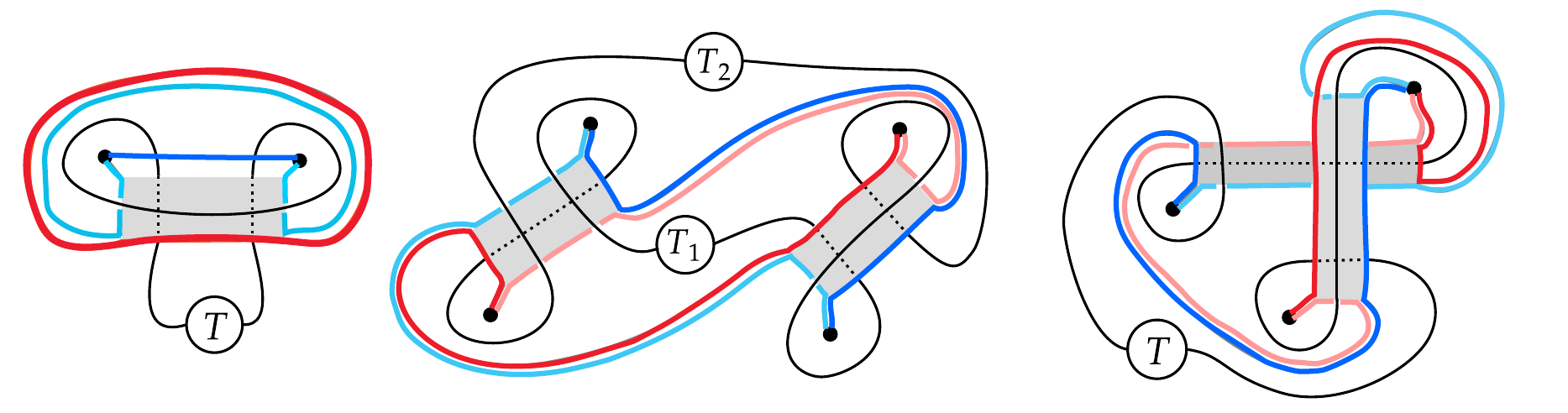}
   \caption{The two boundaries of the annuli appear in blue and red.}
\label{fig:bongleannulibound}
\end{figure}

For all three types, the two boundaries of $A$ together do not separate $\partial H$, so they cannot cobound an annulus on $\partial H$, and hence $A$ is not boundary-parallel. Thus, $A$ is a properly embedded essential annulus, contradicting the hyperbolicity of the bongle.
\end{proof}

In an upcoming paper, we will consider the hyperbolicity of more complicated charm bracelets and present cases in which non-alternating charm bracelets can still be hyperbolic. Additionally, we will show that there exist arbitrarily large nontrivial volume equivalence classes of charm bracelets depending on their self-symmetries.

\section{Volumes of Bongles}



\begin{figure}[htbp]
\includegraphics[scale=0.5]{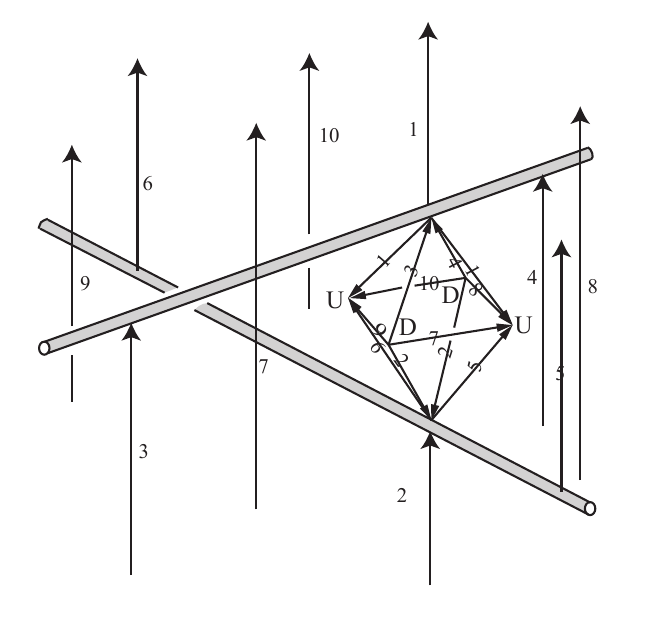}
   \caption{An octahedron at each crossing, with edges labeled to coincide with the edges depicted. Unspecified ends of edges end at the vertex $D$ beneath the projection plane or the vertex $U$ above the projection plane.}
\label{fig:octahedron}
\end{figure}

In \cite{DThurston}, Dylan Thurston introduced a method to decompose a knot or link complement in the 3-sphere into octahedra, each with two ideal vertices and otherwise finite vertices by placing said octahedra at each crossing of a projection of the link. See Figure \ref{fig:octahedron}.
    The vertices labeled $U$ are pulled up to a vertex above the projection sphere and the vertices labeled $D$ are pulled down to a vertex beneath the projection sphere.

\begin{figure}[htbp]
\includegraphics[scale=0.45]{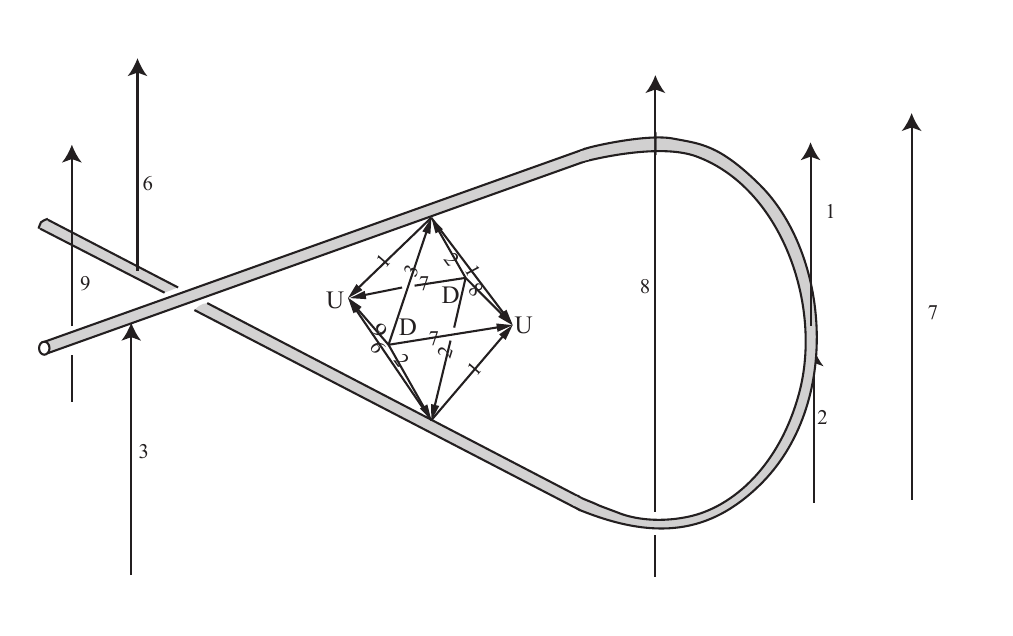}
   \caption{Edges are identified in the case the crossing is in a monogon.}
\label{fig:octahedronmonogon}
\end{figure}

In the case of a bongle, we consider it as a knot in a thickened sphere, with tunnels drilled out for where the stakes go in. Each crossing of the bongle corresponds to a unique monogon in the diagram. Hence, as in Figure \ref{fig:octahedronmonogon}, distinct edge classes as appearing in Figure \ref{fig:octahedron} are identified. In particular, edge classes 7 and 10 coincide, as do edge classes 2 and 4 and edge classes 1 and 5.

Note that the vertices labeled $U$ and $D$, which were finite vertices up above the projection sphere and down below the projection sphere in the case of Thurston's original construction here lie outside the inner and outer spheres in the thickened sphere and hence will be truncated vertices.

In order to remove the stake corresponding to a monogon, we remove an open neighborhood of the vertical edge class corresponding to the center of the monogon, which appears in Figure \ref{fig:octahedronmonogon} as the edge class labeled 8.

\begin{figure}[htbp]
\includegraphics[scale=0.45]{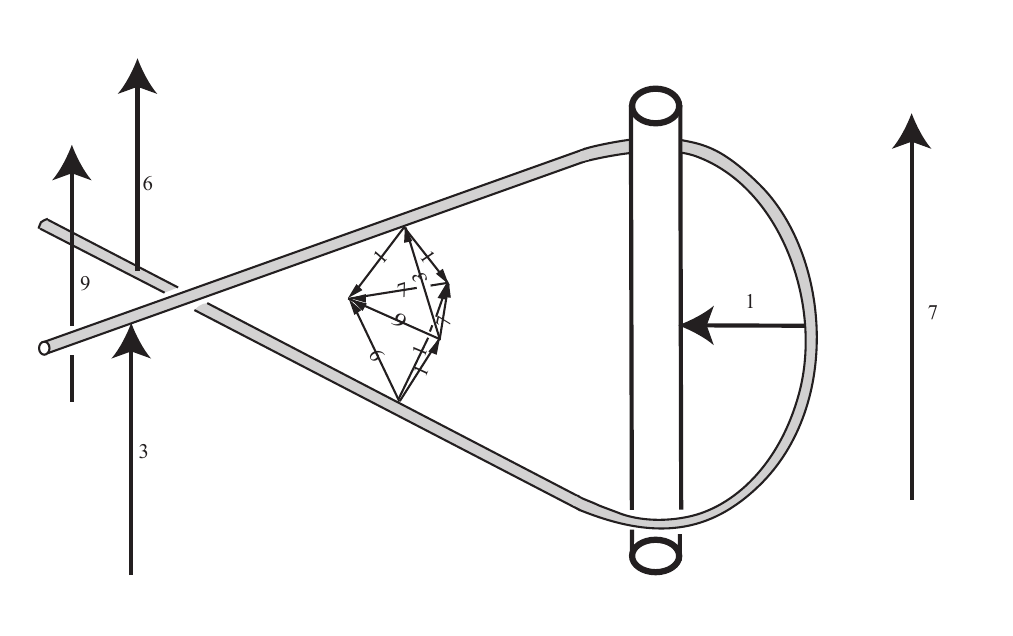}
   \caption{The octahedra become triangular bipyramids when we remove stakes.}
\label{fig:tribipyramid}
\end{figure}

In doing so, we are removing a neighborhood of the edge labeled 8 on the octahedron, and truncating the entire edge. This is topologically equivalent to shrinking that entire edge down to a single truncated vertex. The result is to replace the octahedron by a triangular bipyramid with top and bottom vertices still ideal and the three equatorial vertices all truncated, as in Figure \ref{fig:tribipyramid}. Note that edge classes 1 and 2 are identified by this process.

We now use this construction to prove various results about volumes of bongles.

\begin{theorem}\label{balancedvolume} For fixed $n\ge 4$, all balanced alternating $n$-bongles have the same volume. 
\end{theorem}

\begin{figure}[htbp]
    \centering
    \includegraphics[scale=0.7]{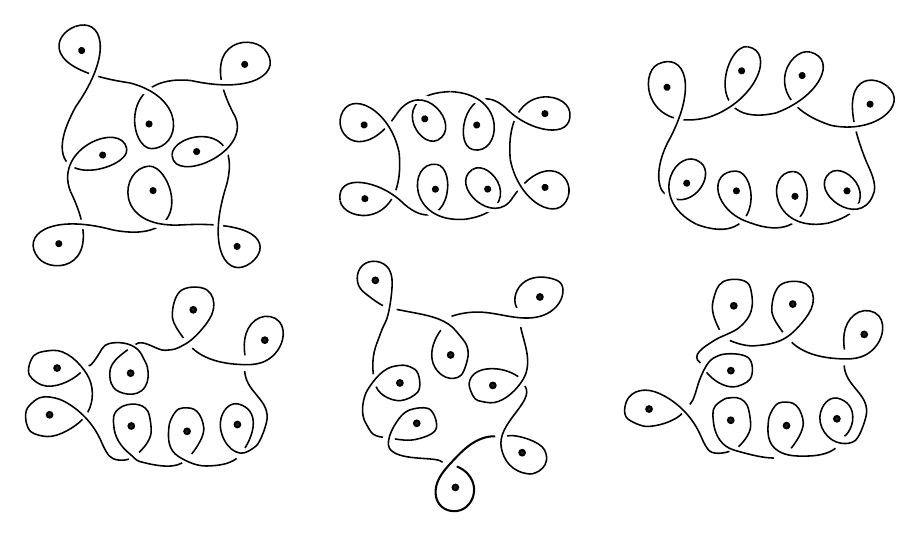}
    \caption{Inequivalent bracelets with the same volume.}
    \label{fig:diffsym}
\end{figure}

For example, the six alternating 8-bongles appearing in Figure \ref{fig:diffsym} all have the same volume.

\begin{proof} Each crossing of a bongle contributes one triangular bipyramid to a decomposition of the complement. For a given bipyramid, the edge class appearing as 1 in Figure \ref{fig:tribipyramid} will appear once each on the triangular bipyramid corresponding to the two adjacent crossings. Hence, each edge class that has an end on an ideal vertex of the bipyramids will have exactly six edges in its edge class. Since the bongle is balanced, there will be a total of $3n/2$ edges in each of the two edge classes corresponding to the central edge and the exterior edge of the bongle. 

If we take all of the triangular bipyramids for the balanced bongle to be isometric with dihedral angles $\pi/3$ for the three edges leading into each ideal vertex and angles $4\pi/3n$ for the equatorial edges, the sum of the angles for each edge class will be exactly $2\pi$.  The ideal triangles are regular, and fit together to satisfy the necessary completeness condition on the cusp. The totally geodesic boundary of the manifold, which is a genus $n-1$ surface, is made up of isometric rhombi, all of the same edge-length, three each from each of the three equatorial vertices of each triangular bipyramid. Each rhombus has two opposite angles of $\pi/3$, and at each such vertex six rhombi meet.  And it has two vertices of angle  $4\pi/3n$, where $3n/2$ rhombi meet. As in \cite{FrigerioPetronio}, this satisfies the necessary conditions to  yield a hyperbolic structure on the manifold of volume given by $n$ times the volume of such a triangular bipyramid. This holds regardless of the order of the innies and outies in the balanced bongle, hence all such balanced $n$-bongles have the same volume.
\end{proof}

\begin{lemma} \label{max3bi} The volume of a 3-bipyramid with the two opposite apices ideal and the other three vertices truncated has supremum 
$5 v_{tet}$, where $v_{tet} \approx 1.01494\dots$ is the volume of an ideal regular tetrahedron.
\end{lemma}
 
Note that we assume the truncation planes are disjoint.  We generalize a similar argument to Proposition 3.3 and Corollary 3.4 of \cite{calderonmayer}.

\begin{proof}
We add a vertical edge from the top ideal vertex to the bottom ideal vertex and then cut the 3-bipyramid into three tetrahedra that share the new edge and that each has two ideal vertices and two truncated vertices. Denote them $T_1, T_2$, and $T_3$. The $k$-th tetrahedron has dihedral angles $A_k, B_k, ..., F_k$ for $k=1,2,3$ with the constraints $$B_k+F_k+D_k=C_k+E_k+D_k=\pi$$ coming from the ideal vertices and an additional constraint $$D_1+D_2+D_3=2\pi$$ coming from the shared central edge.

Note that for any of the three tetrahedra, if $A_k > 0$, we can increase the volume by decreasing $A_k$, so to obtain a maximal volume, we take $A_k = 0$. Although this is no longer a 3-bipyramid with positive distance between the truncation planes, it is a limit of such.

 We use Lagrange multipliers to maximize the volume subject to a loosened version of these constraints. Let
\begin{align*}
g &= \sum_{k = 1}^3 [(B_k +F_k + D_k) +  (C_k + E_k + D_k)]- 6\pi,\\
h &= D_1+D_2+D_3-2\pi.
\end{align*}

The method of Lagrange multipliers with multiple constraints then yields the system of equations given by  
$$\nabla (\vol - \lambda g -  \mu h)= 0,$$
which yields
\begin{align*}
    \frac{\partial \vol}{\partial B_k} &= \frac{\partial \vol}{\partial F_k} = \frac{\partial \vol}{\partial C_k} = \frac{\partial \vol}{\partial E_k} = \lambda\\
    \frac{\partial \vol}{\partial D_k} &= 2\lambda + \mu.
\end{align*}

\medskip

Sch\"afli's Differential formula says that for a generalized hyperbolic tetrahedron $T$, the differential of volume is given by 
$$d\vol(T)= -\frac{1}{2}\sum_{i=1}^6 \ell_i(d\alpha_i),$$ 
where $\alpha_i$ corresponds to the dihedral angle of edge $i$ and $\ell_i$ is the length of edge $i$. Applying this to the above equations yields: 

 $$\ell_{B_1} = \ell_{B_k}=\ell_{F_k} =\ell_{C_k}=\ell_{E_k},$$ 
 $$\ell_{D_k}= 2\ell_{B_1}+ \mu,$$
 
\medskip

\noindent for some constant $\mu$ independent of $k$. Note that although all of these edges have at least one ideal vertex and therefore have infinite length, we can choose finite vertices $\overline{v_1}$ and $\overline{v_2}$ on the vertical edge and measure the distances to them, and take the limit as $\overline{v_1}$ approaches $v_1$ and $\overline{v_2}$ approaches $v_2$.

In particular, all three truncated tetrahedra must be isometric. Therefore, $$D_1=D_2=D_3=2\pi/3$$ and hence, 
 $$B_k=C_k=E_k=F_k=(\pi-D_k)/2 = \pi/6$$ for $k = 1, 2, 3.$ 


 Gluing the three tetrahedra that have these angles together again yields the truncated 3-bipyramid $B_{\infty}$, which has three equatorial edges of angle 0 and all edges with one ideal endpoint having angle $\pi/3.$ Note that for this 3-bipyramid, any pair of truncation planes intersect in a single ideal point. 
 
 \medskip
 
 We now calculate the volume of $B_{\infty}$. 
We first decompose $B_{\infty}$ into two tetrahedra, each with one ideal vertex and three truncated vertices with each pair of truncation surfaces meeting at a point where the edge between the corresponding vertices used to be. After truncation, the resulting polyhedron $P_{\infty}$ has four ideal vertices (the original one and the ones where the truncation planes meet) and three finite vertices. 

From $P_{\infty}$, we can slice off three tetrahedra, all isometric, each with three ideal vertices and one finite vertex, and with angles $\pi/3$ for the edges going into the ideal vertex and angles $\pi/2$ for the remaining two edges going into the finite vertex. Call one of these tetrahedra $R$. Note that after slicing the three copies of $R$ off $P_{\infty}$, we are left with a single ideal regular tetrahedron. 

If we glue six copies of $R$ together around the edge from the ideal vertex to the finite vertex, we obtain a 6-pyramid with totally geodesic 6-gon face on the bottom. Doubling the pyramid across this bottom face yields a 6-bipyramid that decomposes into six ideal regular tetrahedra. Hence, it has volume 6(1.0149416...). Thus, $R$ has volume 1.0149416.../2, and $$ \text{vol} (P_{\infty})= 1.0149416... + \frac{3}{2} (1.0149416...)= \frac{5}{2}(1.0149416\dots).$$ Thus, $$\text{vol}(B_{\infty})= 5(1.0149416\dots)$$.
\end{proof}
 

\begin{theorem}\label{uppervolumebound} For any even $n \geq 2$, the volume of a hyperbolic $n$-bongle is bounded from above by the volume of a balanced hyperbolic $n$-bongle. For $n \geq 2$, every hyperbolic $n$-bongle has volume strictly bounded from above by $5n v_{tet}$. 
    \end{theorem}

    Note that in the case $n$ is even, the volume of $n$ copies of a 3-bipyramid with dihedral angles for the edges going to the ideal vertices of $\pi/3$ and dihedral angles for edges between the truncated vertices of $\frac{4\pi}{3n}$ yields the volume of a balanced $n$-bongle. In the case $n$ is odd, this does not correspond to a volume of an $n$-bongle, but it is still an upper bound on the volume of all hyperbolic $n$-bongles.

\begin{proof} In the case $n = 2$, the hyperbolic balanced 2-bongle is the link exterior $6_2^2$ of volume $4.0597\dots$, and the only other hyperbolic 2-bongle is the Whitehead link exterior of volume $3.6638\dots$, so both  result hold for $n =2$. From now on, we assume $n \geq  3$, and therefore the manifold has totally geodesic boundary of genus at least 2. 

We already have a combinatorial decomposition of the $n$-bongle into $n$ truncated 3-bipyramids. To decompose each 3-bipyramid $B_j$, we add a vertical edge from the top ideal vertex to the bottom ideal vertex and then cut the 3-bipyramid into three tetrahedra that share the new edge. Each tetrahedron has two ideal vertices and two truncated vertices. Denote them $T_{j,1}, T_{j, 2}$, and $T_{j, 3}$. The $(j,k)$-th tetrahedron has dihedral angles $A_{j,k}, B_{j,k}, ..., F_{j,k}$ for $k=1,2,3$ with the constraints $B_{j,k}+F_{j,k}+D_{j,k}=C_{j,k}+E_{j,k}+D_{j,k}=\pi$ coming from the ideal vertices and $D_{j,1}+D_{j,2}+D_{j,3}=2\pi$ coming from the shared central edge.

In particular, these constraints can be loosened to imply:

$$\sum_{j,k}D_{j,k}=2n\pi$$  and $$\sum_{j,k} (B_{j,k}+C_{j,k}+E_{j,k}+F_{j,k})=2n\pi.$$

Note that each edge labelled $A_{j,k}$ corresponds to the outer central edge or the inner central edge of the bongle, so in particular 
$$\sum_{j,k}A_{j,k} = 4\pi.$$ 

Let $$g_1 = \sum_{j,k}D_{j,k}- 2n\pi,$$ $$g_2 = \sum_{j,k} (B_{j,k}+C_{j,k}+E_{j,k}+F_{j,k})- 2n\pi,$$  $$g_3 = \sum_{j,k}A_{j,k} = 4\pi.$$

The method of Lagrange multipliers with multiple constraints then yields the system of equations given by  
$$\nabla (\vol -  \lambda_1 g_{1} -\lambda_2 g_{2} - \lambda_3 g_3)= 0,$$

which yields

$$\frac{\partial \vol}{\partial D_{j,k}} = \lambda_1,$$

$$\frac{\partial \vol}{\partial B_{j,k}} = \frac{\partial \vol}{\partial C_{j,k}} = \frac{\partial \vol}{\partial E_{j,k}} =\frac{\partial \vol}{\partial F_{j,k}} = \lambda_2,$$

$$\frac{\partial \vol}{\partial A_{j,k}} = \lambda_3.$$

As in the previous proof, we apply Sch\"afli's Differential formula to obtain:

$$\ell_{D_{1,1}}=\ell_{D{j,k}},$$
 $$\ell_{B_{1,1}}=\ell_{B_{j,k}} =\ell_{C_{j,k}}=\ell_{E_{j,k}}=\ell_{F_{j,k}},$$ 
 $$\ell_{A_{1,1}}=\ell_{A_{j,k}}.$$

 As in the preceding proof, we must interpret these lengths in terms of limits. 
 
 The equations imply that at this critical point, all $3n$ tetrahedra must be isometric. Moreover, each angle $D_{j,k} = 2\pi/3$ , $A_{j,k} = 4\pi/3n$ and $B_{j,k} = C_{j,k} = E_{j,k} = F_{j,k} = \pi/6$. If we glue the tetrahedra back together to obtain the 3-bipyramids, we obtain the isometric 3-bipyramids we used to construct the balanced $n$-bongles in the proof of Theorem \ref{balancedvolume}.
Thus, when $n$ is even, this critical point corresponds to any of the hyperbolic balanced $n$-bongles.

\medskip

We now prove that this is the global maximum volume over all hyperbolic $n$-bongles. Restrict the domain $D$ of all possible dihedral angles for the $3n$ tetrahedra to lie in the range $(0, \pi)$. Then, we have discovered only one critical point in this range. Note that it is a local maximum value for volume by Lemma 5.3 of \cite{FuterGueritaud}.  We call that volume $V_n^B$. When $n$ is even, it corresponds to a balanced bongle. If we can prove that any hyperbolic $n$-bongle $B_n$ either has a realization of its 3-bipyramids in this same domain, or if it does not have such a realization, it has volume less than $V_n^B$, then this proves that the unique maximum will have volume greater than the volume of that $n$-bongle. 

We can use the hyperbolic structure on the $n$-bongle $B_n$ to straighten the edges corresponding to the 3-bipyramid decomposition relative to the hyperbolic metric. To make this explicit, we begin with one 3-bipyramid and we take one of its ideal vertices and choose it to be the point at $\infty$ in the upper-half-space model of hyperbolic 3-space, denoted $\mathbb{H}^3$. The totally geodesic boundary surface of the $n$-bongle lifts to a disjoint collection of totally geodesic planes in $\mathbb{H}^3$, which appear as hemispheres perpendicular to the horizontal boundary plane.  The three edges coming out of the ideal vertex are realized as three vertical geodesics. These will come down and intersect at right angles three disjoint totally geodesic hemispheres $H_1, H_2$ and $H_3$ that correspond to the truncation planes. 

There will be three unique geodesics that are perpendicular to these planes in pairs. Note that each must have a positive length between the truncation planes, since the boundary surface is totally geodesic. And there will be a unique point in the boundary plane that corresponds to the second ideal vertex of the 3-bipyramid and such that there are three geodesics that share that ideal point and that are each perpendicular to one of $H_1, H_2$, and $H_3$. This unique ideal point is determined by the dihedral angles on the edges ending at this ideal vertex.  This is the geometric realization of the 3-bipyramid. 

There are two facts to note about this geometric realization. First, each of the edges of a 3-bipyramid cannot be isotoped into any of the boundaries of $B_n$, including either the totally geodesic boundary or the torus boundary of the exterior of the knot. This is obvious for the edges that have an ideal endpoint, since they connect the torus boundary to the higher genus boundary. 

For the edges that go from truncated vertex to truncated vertex, they are of two types. Either they correspond to the inner central vertical edge of $B_n$ or they correspond to the outer central vertical edge of $B_n$. These appear as edges labelled 7 and 9 in Figure \ref{fig:octahedronmonogon}. In the case of an $n$-bongle with at least one monogon to the inside and one monogon to the outside, such an edge $e$ bounds a disk $D_1$ together with an arc $\beta$ in the boundary that is punctured once by the knot. Note that $\beta$ intersects a single stake that is in a monogon to the side of the bongle that corresponds to $e$. 

If $e$ were isotopic into the boundary, there would be a disk $D_2$ with boundary $e$ and an arc in the boundary, which after cut and paste could be made to miss $D_1$ on its interior. Hence $A = D_1 \cup D_2$ would be a disk punctured once by the knot. Hence, in the  bongle $B_n$, it becomes an essential annulus, a contradiction to hyperbolicity.

In the case all of the monogons are outies, the explicit hyperbolic structure from \cite{frigerio} implies that all edges are essential and none can be isotoped into the boundary.

Second, note that the geometric 3-bipyramid that we have constructed could be negatively oriented relative to our expectation if any angle is negative.  Or it could be degenerate if an angle is 0 or $\pi$. 

We then move on to the next 3-bipyramid that shares a vertical face with this one and build its geometric realization. In this way, we attempt to build a fundamental domain for the manifold in $\mathbb{H}^3$. If all of the 3-bipyramids are positively oriented, we obtain a valid fundamental domain for the hyperbolic structure on the bongle and the volume of $B_n$ will be exactly the sum of the volumes of the $n$ 3-bipyramids. 

If there are 3-bipyramids that are either degenerate or negatively oriented, meaning this decomposition of $B_n$ does not have a solution in $D$,   it is still the case that the sum of the volumes (including negative values for negatively oriented bipyramids and 0 for degenerate bipyramids)  yields the volume of the manifold.  But, if there is a degenerate or negatively oriented 3-bipyramid,  the volume of $B_n$ is less than $(n-1)(5(1.01494\dots))$ since by Lemma \ref{max3bi}, the  volume of a 3-bipyramid is less than $5(1.01494\dots)$. 

So either there is a realization of the 3-bipyramids for $B_n$ in the domain, in which case $V_n^B$ is greater than the volume of $B_n$ and we are done, or the volume of $B_n$ less than $(n-1)(5(1.01494\dots))$. 

     In fact, we now show $(n-1)(5(1.01494\dots))$ is less than $V_n^B$, so in either case the volume of $B_n$ is less than $V_n^B$, as we wished to show.

The volume $V_n^B$ is $n$ times the volume of the 3-bipyramid with dihedral angles $\pi/3$ going into the two ideal vertices at top and bottom, and dihedral angles $\frac{4\pi}{3n}$ at the three equatorial edges.

Cutting this 3-bipyramid in half equatorially results in two tetrahedra, each with one ideal vertex and three truncated vertices, with three edges of dihedral angle $\pi/3$ going into the ideal vertices and three edges of angle $\frac{2\pi}{3n}$ connecting the truncated vertices. We apply Ushijima's formula for the volume of such a tetrahedron $T_n$ from p. 251-252 of \cite{Ushijima}. We find 
\begin{align*}
\det G &= -6.75 \cos(\frac{2\pi}{3n}),\\
z_1 &= 1 + i \tan\left(2\pi/3n\right),\\
z_2 &= -1 + i \tan\left(2\pi/3n\right),\\
U(z_1, T_n) &= \frac{3}{2}\left[\Li_2(e^{\frac{2\pi i}{3}} e^{\frac{4\pi i}{3n}}(1 + i \tan\left(2\pi/3n\right)) - \Li_2(-e^{\frac{\pi i}{3}} e^{\frac{4\pi i}{3n}}(1 + i \tan\left(2\pi/3n\right)\right],\\
U(z_2, T_n) &= \frac{3}{2}\left[ \Li_2(e^{\frac{2\pi i}{3}} e^{\frac{4\pi i}{3n}}(-1 + i \tan\left(2\pi/3n\right)) - \Li_2(-e^{\frac{\pi i}{3}} e^{\frac{4\pi i}{3n}}(-1 + i \tan\left(2\pi/3n\right))\right],\\
\Vol(T_n) &= \frac{1}{2} \Im [U(z_1, T) - U(z_2, T)].
\end{align*}

Then, we want to show that $$V_n^B > (n-1) 5.0747\dots.$$ But, $V_n^B = 2n \Vol(T_n)$, so we want to show that $$\frac{2n}{n-1}\Vol(T_n) > 5.0747\dots$$ for all $n \geq 3$.
Note that as $n$ approaches $\infty$, the left hand side approaches the right. 

As can be seen using Mathematica, the  function $\frac{n}{n-1}\Im [U(z_1, T_n)]$ is decreasing for all $n \geq 12$. The function
$\frac{n}{n-1}\Im [U(z_2, T_n)]$ is decreasing for all $n \geq 1$. Thus, at least for $n \geq 12$, the function on the left is decreasing toward 5.0747\dots. Hence, it must be greater than 5.0747\dots.

For $3 \leq n \leq 11$, we can calculate the volume of $T_n$ and see that $$\frac{2n}{n-1}\Vol(T_n) > 5.0747\dots.$$ So, this inequality holds for all $n \geq 3$.
\end{proof}

Note that any alternating $n$-bongle with $k$ innies and $n-k$ outies has $2k + n-k = k+ n$ edges in the edge class corresponding to the central edge and has $k + 2 (n-k) = 2n -k$ edges in the edge class corresponding to the exterior edge. We may make a choice of which is the inside and which is the outside of the bongle so that $k\le n/2$. 

\begin{conjecture} \label{volumeorder} For hyperbolic $n$-bongles with a fixed value of $n \ge 3$, as $k$ increases toward $n/2$, volume increases.
\end{conjecture}

\begin{conjecture}\label{limitvolume}
Let $(v_3, v_4, v_5, \dots)$ be a sequence such that $v_n$ is the volume of some hyperbolic $n$-bongle. Then $$\lim_{n \to \infty} v_n/n = 5.0747 \dots.$$
\end{conjecture}

Note that this conjecture is known to be true in the case all of the $n$-bongles in the sequence are alternating bongles that are all outies. This was proved in \cite{frigerio} using explicit angles to represent said bongles. However, there, a calculation error lead to a statement with a different value. Thanks to R. Frigerio for help in determining the correct value. 
     Similarly, this is known in the case the bongles in the sequence are either all balanced or a mixture of balanced and all outie bongles. 

In order to prove the full conjecture, by  Frigerio's result and Theorem \ref{uppervolumebound}, it would be enough to show that an alternating all-outie $n$-bongle has least volume among all alternating $n$-bongles for $n \geq 3$. This would follow from the stronger Conjecture \ref{volumeorder}.






\newpage

\section{Volume Table}

Note that we do not show the other balanced 4-bongle and 6-bongles since they have the same volume as the balanced bongle shown.

\begin{figure}[htbp]
    \centering
    \includegraphics[scale=.3]{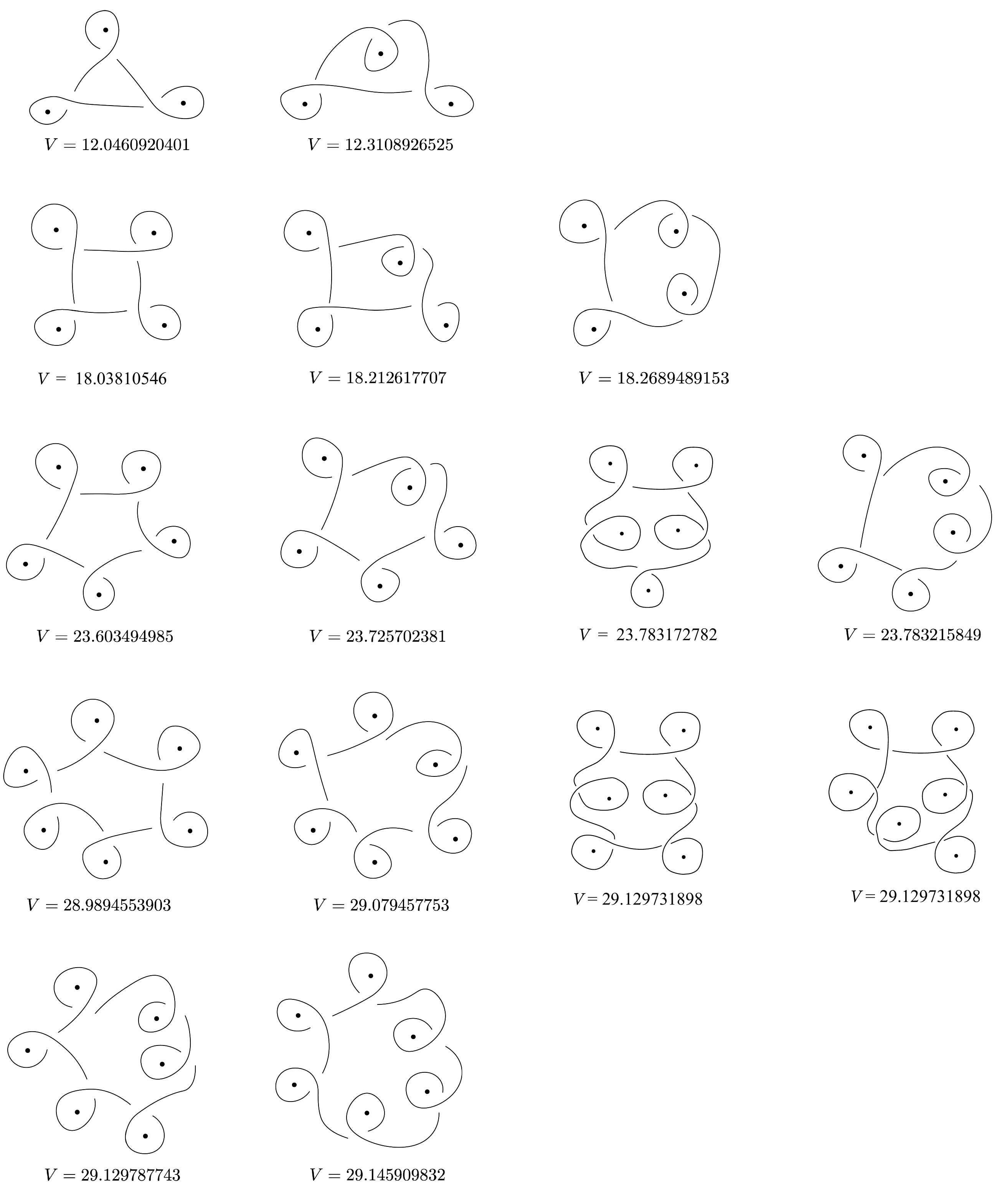}
    \caption{Volumes of alternating 3-bongles through 6-bongles.}
    \label{fig:bongle chart}
\end{figure}

\bibliographystyle{plain}
\bibliography{bib}

\end{document}